\documentclass{amsart}
\usepackage{amsbsy,amssymb}
\usepackage{epsfig}
\usepackage{graphicx}

\title{Unbiased Estimators for Entropy and Class Number}

\author{Stephen Montgomery-Smith}
\address{Math Dept, University of Missouri, Columbia, MO 65211, U.S.A.}
\email{stephen@missouri.edu}
\author{Thomas Sch\"urmann}
\address{J\"ulich Supercomputing Centre, J\"ulich Research Centre, 52425 J\"ulich, Germany.}
\email{thomas.schuermann@live.de}

\begin{document}

\begin{abstract}We introduce unbiased estimators for the Shannon entropy 
and the class number, in the situation that
we are able to take sequences of independent samples of arbitrary 
length.\end{abstract}

\maketitle

\subsection*{Introduction}

This paper supposes that we may pick a sequence of arbitrary length
of independent samples 
$w_1,w_2,\dots$ from
an infinite population.  
Each sample belongs to one of $M$ classes $C_1,C_2,\dots,C_M$, 
and the probability that a sample belongs to
class $C_i$ is $p_i$.  So these probabilities
satisfy the constraints $0\le p_i\le 1$ and $\sum_{i=1}^M\,p_i =1$.

The goal of this paper is to present methods to 
estimate the Shannon entropy $H=-\sum_{i=1}^M\,p_i\log(p_i)$
(see \cite{shann}).
An obvious method is to take a sample of size $n$, and compute the 
estimators $\hat p_i = k_i/n$, where $k_i$ are the number
of samples from the class $C_i$.  
However this is known to systematically underestimate the entropy, 
and it can be significantly
biased \cite{grass88}\cite{harris}\cite{herzel}\cite{miller}.  
Recently there have been more advanced estimators for the entropy
which have smaller bias 
\cite{grass88}\cite{grass03}\cite{harris}\cite{herzel}\cite{miller}\cite{schuer}. 

In this paper we introduce new entropy estimators that have bias identically 
zero.  
We also introduce an unbiased estimator for the class number $M$, a problem
of interest to ecologists (see for example the review article \cite{chao05}).
The disadvantage of all our methods is that there is no
\emph{a priori} estimate of the sample size.  
For this reason, we postpone rigorous
analysis of variance and other measures of confidence until it 
becomes clear that these estimators are of more than theoretical value.

We will use the following power series.  Define the harmonic number by
$h_n = \sum_{k=1}^n 1/k$, $h_0 = 0$.  Then for $|x|<1$
\begin{eqnarray*} 
   \frac1{(1-x)^2} = \displaystyle \sum_{k=1}^\infty k x^{k-1} ,
\\
   \log(1-x) = - \displaystyle \sum_{k=1}^\infty \frac{x^k}k ,
\\
   \frac{\log(1-x)}{(1-x)} = - \displaystyle \sum_{k=1}^\infty h_k x^k .
\end{eqnarray*}

\subsection*{First Estimator for Entropy}

For each $1 \le i \le M$, let $N_i$ denote the smallest $k \ge 1$ 
for which $w_k \in C_i$.  Then
$$ \hat H_1 = \sum_{i=1}^M \frac{I_{N_i\ge 2}}{N_i-1} $$
is an unbiased estimator for the entropy.  The proof is straightforward.  
The marginal distribution of $N_i$ satisfies the
geometric distribution $\Pr(N_i = k) = p_i (1-p_i)^{k-1}$.  Thus
$$ E\left(\frac{I_{N_i\ge 2}}{N_i-1}\right) = 
\sum_{k=2}^\infty \frac{p_i(1-p_i)^{k-1}}{k-1} = - p_i \log(p_i) .$$

Depending upon the applications, a possible disadvantage of this estimator 
is that complete knowledge of all possible classes
needs to be known in advance.

\subsection*{Second Estimator for Entropy}

Let $N$ denote the smallest $k \ge 1$ such that $w_1$ and
$w_{k+1}$ belong to the same class.  
Then $\hat H_2 = h_{N-1}$ is an unbiased estimator for the entropy.

This follows, since conditional upon $w_1 \in C_i$, 
the distribution of $N$ satisfies the geometric distribution
$\Pr(N = k | w_1 \in C_i) = p_i (1-p_i)^{k-1}$.  Thus
$$ E(h_{N-1} | w_1 \in C_i) = \sum_{k=1}^\infty h_{k-1} p_i (1-p_i)^{k-1} 
   = - \log(p_i) ,$$
and hence
$$ E(h_{N-1}) 
   = \sum_{i=1}^M E(h_{N-1} | w_1 \in C_i) \Pr(w_1 \in C_i) 
   = - \sum_{i=1}^M p_i \log(p_i) .$$

While we don't wish to focus on analysis of the variance, it is certainly 
clear that this single estimator by itself will have
unusable confidence limits.  While the variance can be reduced by taking the 
mean of $n$ of these estimators, we
propose the following version.
For each $1 \le j \le n$, let 
$N^{(j)}$ be the smallest $k \ge 1$ such that $w_j$ and 
$w_{k+j}$ belong to the same class.  Then
the unbiased estimator is
$$ \hat H_3 = \frac1n \sum_{j=1}^n h_{N^{(j)}-1} .$$

\subsection*{Estimator for Class Number}

This is very similar to the second estimator for entropy.  
Define $N$ and $N^{(j)}$ as in the previous section.  
Then $\hat M_1 = N$ is an unbiased estimator for the class number $M$.
The proof is almost identical to that provided in the previous section.

We may also produce an unbiased estimator
$$ \hat M_2 = \frac1n \sum_{j=1}^n N^{(j)} .$$
This last quantity can also be considered as a corrector to the
naive estimator $\hat M_3$, which is defined as the number of classes
observed
in the first $n$ samples, or alternatively, as
the cardinality of the
set $A$, where
$A$ is the set of $1 \le i \le M$ such that there exists
$1 \le k \le n$ for which $w_k \in C_i$.

After picking $n$ samples, then continue picking samples until every
class observed in the first $n$ samples is observed at least once more.
For each $i \in A$, let $F_i$ denote the smallest $k \ge 1$ such that 
$w_k \in C_i$,
and let $L_i$ denote the smallest $k \ge 1$ such that $w_{k+n} \in C_i$.  
Then
$$ \hat M_2 = \hat M_3 + \frac1n \sum_{i\in A} (L_i - F_i) .$$

If one doesn't wish to record the order in which the samples are obtained, one can simply compute the expected value of this quantity over all possible
rearrangements of obtaining this data, and derive the unbiased estimator
$$ \hat M_4 = \hat M_3 + \frac1n \sum_{i\in A}
   \left(\frac{m+1}{s_i+1}-\frac{n+1}{r_i+1} \right) ,$$
where $r_i$ is the number of times the $i$th class appears in the first 
sample
of size $n$, $m$ is the size of the subsequent sample, and $s_i$ is the
number of times the $i$th class appears in the subsequent sample.

\end{document}